\documentclass[12pt,twoside,spanish,a4paper]{article}

\usepackage{graphicx}
\usepackage{epstopdf}
\usepackage[T1]{fontenc}
\newcommand{\medpar}{\par\medskip}
\newcommand{\Eqref}[1]{(\ref{#1})}

\title{The Pearcey integral in the highly oscillatory region}
\author{\\Jos\'{e} L. L\'{o}pez$^1$ and Pedro Pagola$^2$\\\\
\small{$^1$
\textsf{\textit{Dpto. de Ingenier\'{\i}a Matem\'{a}tica e Inform\'{a}tica and INAMAT, Universidad P\'{u}blica de Navarra}}}\\
\small{ \textsf{\textit{ e-mail: jl.lopez@unavarra.es}}}
\\\
\small{$^2$
\textsf{\textit{Dpto. de Ingenier\'{\i}a Matem\'{a}tica e Inform\'{a}tica, Universidad P\'{u}blica de Navarra}}}\\
\small{ \textsf{\textit{ e-mail: pedro.pagola@unavarra.es}}}}
\date{}

\begin{document}
\normalsize \maketitle

\begin{abstract}
We consider the Pearcey integral $P(x,y)$ for large values of $\vert y\vert$ and bounded values of $\vert x\vert$. The integrand of the Pearcey integral oscillates wildly in this region and the asymptotic saddle point analysis is complicated. Then we consider here the modified saddle point method introduced in [Lopez, P\'erez and Pagola, 2009]. With this method, the analysis is simpler and it is possible to derive a complete asymptotic expansion of $P(x,y)$ for large $\vert y\vert$. The asymptotic analysis requires the study of three different regions for $\arg y$ separately. In the three regions, the expansion is given in terms of inverse powers of $y^{2/3}$ and the coefficients are elementary functions of $x$. The accuracy of the approximation is illustrated with some numerical experiments.
 \\\\
\noindent \textsf{2010 AMS \textit{Mathematics Subject
Classification:} 33E20; 41A60.
} \\\\
\noindent  \textsf{Keywords \& Phrases:} Pearcy integral. Asymptotic expansions. Modified saddle point method.
\\\\
\end{abstract}

\section{Introduction}

The mathematical models of many short wavelength phenomena, specially wave propagation and optical diffraction, contain, as a basic ingredient, oscillatory integrals with several nearly coincident stationary phase or saddle points.
The uniform approximation of those integrals can be expressed in terms of certain canonical integrals and their derivatives \cite{connor}, \cite{ursell}. The importance of these canonical diffraction integrals is stressed in \cite{paris} by means of a very appropriate sentence: {\it The role played by these canonical diffraction integrals in the analysis of caustic wave fields is analogous to that played by complex exponentials in plane wave theory}.

Apart from their mathematical importance in the uniform asymptotic approximation of oscillatory integrals \cite{olde}, the canonical diffraction integrals have physical applications in the description of surface gravity waves \cite{kelvin}, \cite{ursel}, bifurcation sets, optics, quantum mechanics and acoustics (see \cite[Sec. 36.14]{nist} and references there in).

In \cite[Chap. 36]{nist} we can find a large amount of information about this integrals. First of all, they are classified according to the number of free independent parameters that describe the type of singularities arising in catastrophe theory, that also corresponds to the number of saddle points of the integral. The simplest integral with only one free parameter, that corresponds to the fold catastrophe, involves two coalescing stationary points: the well-known integral representation of the Airy function. The second one, depending on two free parameters corresponds to the cusp catastrophe and involves three coalescing stationary points. The canonical form of the oscillatory integral describing the cusp diffraction catastrophe is given by the Cusp catastrophe or Pearcy integral \cite[p.777, eq. 36.2.14]{nist}:
\begin{equation}\label{pea}
\bar P(x,y):=\int_{-\infty}^\infty e^{i(t^4 +xt^2+yt)}dt.
\end{equation}
This integral was first evaluated numerically by using quadrature formulas in \cite{pearcey} in the context of the investigation of the electromagnetic field near a cusp. 
The third integral of the hierarchy is the Swallowtail integral that depends on three free parameters and involves four coalescing stationary points. Apart from the classification of this family of integrals, in \cite[Chap. 36]{nist} we can find many properties such as symmetries, ilustrarive pictures, bifurcation sets, scaling relations, zeros, convergent series expansions, differential equations and leading-order asymptotic approximations among others. But we cannot find many details about asymptotic expansions.

The three first canonical integrals: Airy function, Pearcy integral and Swallowtail integral are the most important ones in applications. The first one is well-known and deeply investigated. In this paper we focus our attention in the second one and will consider the third one in a further research.
The integral \Eqref{pea} exists only for $0\le\arg{x}\le \pi$ and real $y$. As it is indicated in \cite{paris}, after a rotation of the integration path through an angle of $\pi/8$ that removes the rapidly oscillatory term $e^{it^4}$, the Pearcy integral may be written in the form $\bar P(x,y)=2e^{i\pi/8}P(xe^{-i\pi/4},ye^{i\pi/8})$, with
\begin{equation}\label{pearcy}
P(x,y):=\int_0^\infty e^{-t^4 -xt^2}\cos(yt)dt.
\end{equation}
This integral is absolutely convergent for all complex values of $x$ and $y$ and represents the analytic continuation of the Pearcy integral $\bar P(x,y)$ to all complex values of $x$ and $y$ \cite{paris}. Therefore, it is more convenient to work with the representation \Eqref{pearcy} of the Pearcy integral.

We can find in the literature several asymptotic expansions of $P(x,y)$ in different regions of $(x,y)$.
In \cite{kaminski} we can find an asymptotic expansion of the Pearcey integral when $(x,y)$ are near the caustic $8x^3-27y^2=0$ that remains valid as $\vert x\vert\to\infty$. The expansion is given in terms of Airy functions and its derivatives and the coefficients are computed recursively. We refer the reader to \cite{kaminski} for further details.

An exhaustive  asymptotic analysis of this integral can be found in \cite{paris}. In particular, a complete asymptotic expansion for large $\vert x\vert$ is given in \cite{paris} by using asymptotic techniques for integrals applied to the integral \Eqref{pearcy}.
The integral $P(x,y)$ is also analyzed in \cite{paris} for large $\vert y\vert$. But the analysis derived from the standard saddle point method is cumbersome and only the first order term of the asymptotic expansion is given. Then, we are interested here in the derivation of a complete asymptotic expansion of $P(x,y)$ for large values of $\vert y\vert$.

In the following section we analyze the saddle point features of the Pearcey integral for large $\vert y\vert$. In Section 3 we derive a complete asymptotic expansion of $P(x,y)$ for large $\vert y\vert$. Section 4 contains some numerical experiments and a few remarks. Through the paper we use the principal argument $\arg z\in(-\pi,\pi]$ for any complex number $z$.

\section{The saddle point analysis of the Pearcey integral}

Because $P(x,y)=P(x,-y)$, without loss of generality, we may restrict ourselves to the half plane $\Re y\ge 0$. We write the Pearcey integral $P(x,y)$ in the form
$$
P(x,y)=\frac{1}{2}\int_{-\infty}^{\infty} e^{-u^4-xu^2+iyu}du.
$$
Define $\theta:=\arg y$ (it is restricted to $\vert\theta\vert\le\pi/2$). After the change of variable $u=ty^{1/3}=t\vert y\vert^{1/3}e^{i\theta/3}$ we find that
\begin{equation}\label{pxy}
P(x,y)=\displaystyle{\frac{y^{1/3}}{2}\int_{-\infty e^{-i\theta/3}}^{\infty e^{-i\theta/3}} e^{\vert y\vert^{4/3}f(t)-xy^{2/3}t^2}dt, }
\end{equation}
with phase function $f(t):=e^{4i\theta/3}(it-t^4)$. This phase function has three saddle points:
$$
t_0=-\frac{i}{4^{1/3}}, \hskip 3cm t_1=\frac{e^{i\pi/6}}{4^{1/3}}, \hskip 3cm t_2=\frac{e^{5i\pi/6}}{4^{1/3}}.
$$
From the steepest descent method \cite[Chap. 2]{wong}, \cite{saddle}, we know that the asymptotically relevant saddle points are those ones for which the integration path $C:=\lbrace r e^{-i\theta/3}$; $-\infty<r<\infty\rbrace$ in \Eqref{pxy} can be deformed to a steepest descent path (or union of steepest descent paths) that contains a saddle point (or several saddle points). This is only possible for the saddle points $t_1$ and $t_2$. Then, following the idea introduced in \cite{saddle}, we rewrite the phase function $f(t)$ in the form of a Taylor polynomial at the saddle points $t_1$ and $t_2$:
\begin{equation}\label{efe}
\begin{array}{cl}
f(t)= & \displaystyle{\frac{3e^{i(4\theta+2\pi)/3}}{4^{4/3}}-\frac{3e^{i(4\theta+\pi)/3}}{2^{1/3}}(t-t_1)^2-2^{4/3}e^{i(4\theta/3+\pi/6)}(t-t_1)^3-e^{i4\theta/3}(t-t_1)^4, }\\\\
f(t)= & \displaystyle{\frac{3e^{i(4\theta-2\pi)/3}}{4^{4/3}}-\frac{3e^{i(4\theta-\pi)/3}}{2^{1/3}}(t-t_2)^2+2^{4/3}e^{i(4\theta/3-\pi/6)}(t-t_2)^3-e^{i4\theta/3}(t-t_2)^4.}
\end{array}
\end{equation}
From \cite{saddle}, we know that it is not necessary to compute the steepest descent paths of $f(t)$ at $t_1$ and $t_2$, but only the steepest descent paths of the quadratic part of $f(t)$ at those points, that are simpler: they are nothing but straight lines. The steepest descent path of the quadratic part of $f(t)$ at $t_1$ is the straight $C_1:=\lbrace t_1+r e^{-i(\pi+4\theta)/6}$; $-\infty<r<\infty\rbrace$. The steepest descent path of the quadratic part of $f(t)$ at $t_2$ is the straight $C_2:=\lbrace t_2+r e^{i(\pi-4\theta)/6}$; $-\infty<r<\infty\rbrace$.

From \cite{saddle} we also know that it is not necessary to deform the path $C$ into one of these paths or the union of them. It is enough to deform $C$ into a new integration path $\Gamma$ that contains a portion of $C_1$ and/or $C_2$ that includes $t_1$ and/or $t_2$ respectively. The new integration path $\Gamma$ depends on the value of $\theta$. We distinguish the three cases that we specify in the following pictures:

\begin{figure}
\begin{center}
\includegraphics[width=0.65\textwidth]{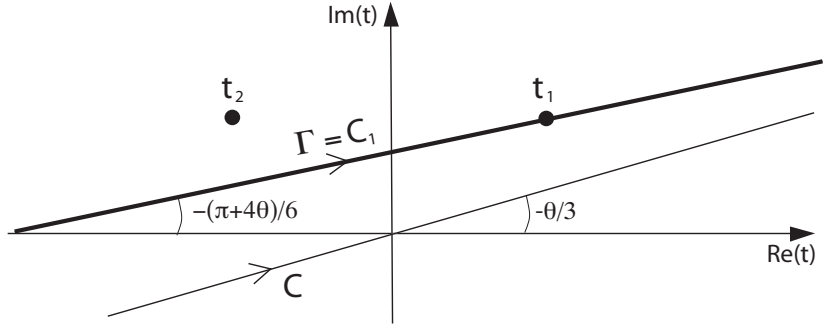}
\end{center}
\caption{\small Case 1. When $\displaystyle{-\frac{\pi}{2}\le\theta< -\frac{\pi}{8}}$, the new integration path $\Gamma$ is just the steepest descent path $C_1$ of the quadratic part of $f(t)$ at $t=t_1$.}
\label{fig:path2}
\end{figure}

\begin{figure}
\begin{center}
\includegraphics[width=0.65\textwidth]{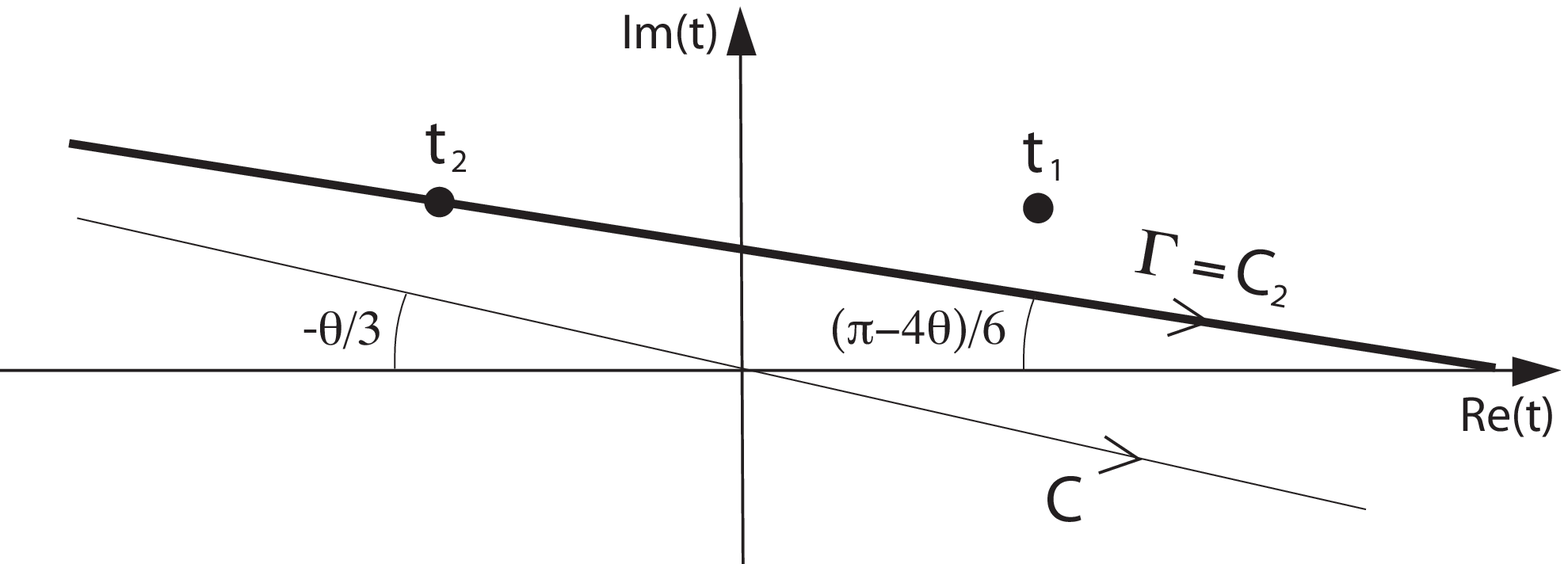}
\end{center}
\caption{\small Case 2. When $\displaystyle{\frac{\pi}{8}<\theta\le \frac{\pi}{2}}$, the new integration path $\Gamma$ is just the steepest descent path $C_2$ of the quadratic part of $f(t)$ at $t=t_2$.}
\label{fig:path3}
\end{figure}

\begin{figure}
\begin{center}
\includegraphics[width=0.65\textwidth]{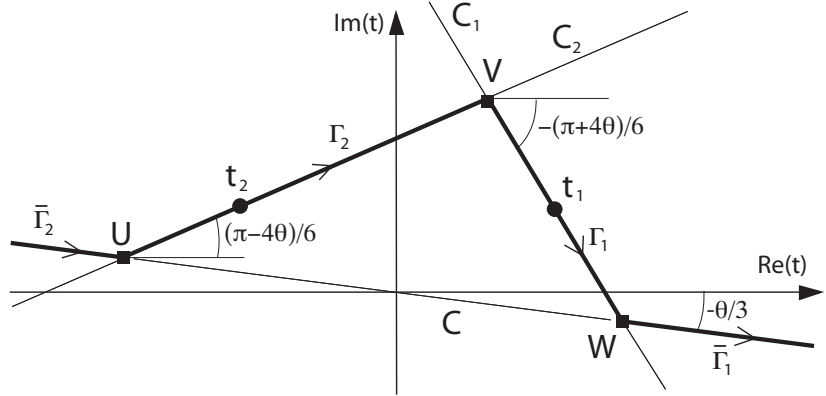}
\end{center}
\caption{\small Case 3. When $\displaystyle{-\frac{\pi}{8}\le\theta\le\frac{\pi}{8}}$, the new integration path $\Gamma$ is the union of the two finite segments $\Gamma_1\subset C_1$ and $\Gamma_2\subset C_2$ and the two infinite segments $\bar\Gamma_1\subset C$ and $\bar\Gamma_2\subset C$.}
\label{fig:path1}
\end{figure}

In any of the three cases, the deformation of the original integration path $C$ to the new path specified in Figures 1-3 can be justified by using the Cauchy's residua theorem. In case 3 it is obvious. In case 1 it follows from the fact that $\vert\arg(e^{4i\theta/3}t^4)\vert<\frac{\pi}{2}$ for $-\frac{\pi+4\theta}{6}\le\arg t\le -\frac{\theta}{3}$ when $-\frac{\pi}{2}\le\theta<-\frac{\pi}{8}$. In case 2 it follows from the fact that $\vert\arg(e^{4i\theta/3}t^4)\vert<\frac{\pi}{2}$ for $-\frac{\theta}{3}\le\arg t\le\frac{\pi-4\theta}{6}$ when $\frac{\pi}{8}<\theta\le\frac{\pi}{2}$.

In case 3, to apply the method introduced in \cite{saddle}, we need to show that the contribution of the integral on the paths $\bar\Gamma_1$ and $\bar\Gamma_2$ is negligible. In fact, the maximum of the real part of $f(t)$ over the paths $\bar\Gamma_1$ and $\bar\Gamma_2$ is attained at the points $\displaystyle W=t_1+2^{-2/3}e^{-i\frac{\pi+4\theta}{6}}$ and $\displaystyle U=t_2-2^{-2/3}e^{i\frac{\pi-4\theta}{6}}$ respectively (see Fig. 3). We have that $\Re f(U)$, $\Re f(W)\le -1.38077...$, that is the value of $\Re f(U)$ for $\theta=\pi/8$ or the value of $\Re f(W)$ for $\theta=-\pi/8$. Therefore we have that,

\begin{equation}\label{negli}
\displaystyle{\int_{\bar\Gamma_1\bigcup \bar\Gamma_2} e^{\vert y\vert^{4/3}f(t)-xy^{2/3}t^2}dt=\mathcal{O}(e^{-1.38077\vert y\vert^{4/3}})},\hskip 3cm \vert y\vert\to\infty.
\end{equation}
\section{A complete asymptotic expansion of the Pearcey integral}

We will see that, in case 3, the integral \Eqref{negli} is negligible (exponentially small) compared with the integrals on $\Gamma_1$ and $\Gamma_2$. Therefore, apart from the exponentially negligible term \Eqref{negli}, the integral \Eqref{pxy} over the path $C$ equals the integral over the path $\Gamma_1\bigcup\Gamma_2$ in case 3, the integral over the path $C_1$ in case 1 or the integral over the path $C_2$ in case 2. From \cite{saddle} we know that, in case 3, the integral over every path $\Gamma_k$, $k=1,2$, can be replaced by an integral over the respective path $C_k$ plus exponentially small terms. Then, in this section, we derive an asymptotic expansion of the integrals on $C_1$ and $C_2$ in any of the three cases. We define
$$
P_k(x,y):=\frac{y^{1/3}}{2}\int_{C_k} e^{\vert y\vert^{4/3}f(t)-xy^{2/3}t^2}dt; \hskip 3cm
k=1,2.
$$
To compute $P_1(x,y)$ in cases 1 and 3 we introduce the change of variable $t=t_1+u\,y^{-2/3}e^{-i\pi/6}$ and consider the form given in the first line of \Eqref{efe} for $f(t)$:
\begin{equation}\label{p1}
P_1(x,y)=\frac{1}{2y^{1/3}}e^{3\cdot 4^{-4/3}e^{2i\pi/3}y^{4/3}-4^{-2/3}xy^{2/3}e^{i\pi/3}-i\pi/6}\int_{u_-}^{u_+} e^{-3\cdot 2^{-1/3}u^2-2^{1/3}xu}e^{h_1(u,x,y)}du,
\end{equation}
with
$$
h_1(u,x,y):=\frac{e^{2i\pi/3}u^2(x+2^{4/3}u)}{y^{2/3}}+\frac{e^{i\pi/3}u^4}{y^{4/3}},
$$
$(u_-,u_+)=(-\infty,\infty)$ in case 1 and $\displaystyle (u_-,u_+)=\left(-(2\vert y\vert^2)^{1/3}\cos\left(\frac{\pi+2\theta}{3}\right),\left(\frac{\vert y\vert}{2}\right)^{2/3}\right)$ in case 3.

To compute $P_2(x,y)$ in cases 2 and 3 we introduce the change of variable $t=t_2+uy^{-2/3}e^{i\pi/6}$ and write $f(t)$ in the form given in the second line of \Eqref{efe}:
\begin{equation}\label{p2}
P_2(x,y)=\frac{1}{2y^{1/3}}e^{3\cdot 4^{-4/3}e^{-2i\pi/3}y^{4/3}-4^{-2/3}xy^{2/3}e^{-i\pi/3}+i\pi/6}\int_{u_-}^{u_+} e^{-3\cdot 2^{-1/3}u^2+2^{1/3}xu}e^{h_2(u,x,y)}du,
\end{equation}
with
$$
h_2(u,x,y):=\frac{e^{-2i\pi/3}u^2(x-2^{4/3}u)}{y^{2/3}}+\frac{e^{-i\pi/3}u^4}{y^{4/3}},
$$
$(u_-,u_+)=(-\infty,\infty)$ in case 2 and $(u_-,u_+)=\displaystyle{\left(-\left(\frac{\vert y\vert}{2}\right)^{2/3},(2\vert y\vert^2)^{1/3}\cos\left(\frac{\pi-2\theta}{3}\right)\right)}$ in case 3.

In order to find an asymptotic expansion of $P_1(x,y)$ and $P_2(x,y)$ for large $\vert y\vert$, it is enough to expand the exponentials $e^{h_1(u,x,y)}$ in \Eqref{p1} and $e^{h_2(u,x,y)}$ in \Eqref{p2} in inverse powers of $y^{2/3}$ \cite{saddle}. To this end we write
$$
e^{h_1(u,x,y)}=\sum_{n=0}^\infty\frac{[h_1(u,x,y)]^n}{n!},\hskip 3cm
e^{h_2(u,x,y)}=\sum_{n=0}^\infty\frac{[h_2(u,x,y)]^n}{n!}.
$$
Applying the binomial Newton formula to $[h_{1}]^n$ and $[h_{2}]^n$ and rearranging the summation indexes we find that
\begin{equation}\label{expah}
e^{h_1(u,x,y)}=\sum_{n=0}^\infty\frac{\bar A_n(x,u)}{y^{2n/3}},\hskip 2cm
e^{h_2(u,x,y)}=\sum_{n=0}^\infty\frac{\bar B_n(x,u)}{y^{2n/3}},
\end{equation}
with
$$
\bar A_n(x,u):=e^{\frac{-i\,n\pi}{3}}\hspace{-0.2cm}\sum_{m=\lfloor\frac{n+1}{2}\rfloor}^n\sum_{k=0}^{2m-n}a_{n,m,k}(x)u^{2m+n-k},
$$
$$
\bar B_n(x,u):=e^{\frac{i\,n\pi}{3}}\hspace{-0.2cm}\sum_{m=\lfloor\frac{n+1}{2}\rfloor}^n\sum_{k=0}^{2m-n}a_{n,m,k}(x)(-u)^{2m+n-k},
$$
\begin{equation}\label{aa}
a_{n,m,k}(x):=\frac{x^k\,2^{4(2m-n-k)/3}\,(-1)^{m}}{k!(2m-n-k)!(n-m)!}.
\end{equation}
Introducing the expansions \Eqref{expah} in \Eqref{p1} and \Eqref{p2} and interchanging sum and integral we find, for $k=1,2$ \cite{saddle}:

\begin{equation}\label{puno}
P_k(x,y)\sim \frac{\sqrt{\pi/3}}{2^{5/6}y^{1/3}}{\rm Exp}\left[{3\frac{y^{4/3}}{4^{4/3}}e^{-(-1)^{k}2i\frac{\pi}{3}}-x\frac{y^{2/3}}{4^{2/3}}e^{-(-1)^{k}i\frac{\pi}{3}}+\frac{x^2}{6}}\right]
\sum _{n=0}^{\infty} e^{ (-1)^k(2 n+1)\frac{ i \pi }{6}}\,\frac{A_n(x)}{y^{2 n/3}},
\end{equation}
with
\begin{equation}\label{pdos}
A_n(x):=\hspace{-0.4cm}\sum _{m=\left[\frac{n+1}{2}\right]}^n \sum _{k=0}^{2m-n} (-1)^{n+k}a_{n,m,k}(x)\,c_{2m+n-k}(x)
\end{equation}
and
\begin{equation}\label{ce}
c_n(x):=\frac{1}{2^{1/6}}\sqrt{\frac{3}{\pi}}e^{-x^2/6}\int_{-\infty}^{\infty} e^{-3\cdot 2^{-1/3}u^2+2^{1/3}xu}u^{n}du=\frac{x^n}{3^n\,2^{n/3}}\sum_{k=0}^{\left[\frac{n}{2}\right]}\left(\frac{3}{2\,x^2}\right)^k\frac{n!}{k!(n-2k)!}.
\end{equation}
The coefficients $c_n(x)$ may be also computed recursively in the form:
$$
\begin{array}{cl}
& c_0(x)=1; \hskip 2cm c_1(x)=\displaystyle{\frac{x}{3\cdot 2^{1/3}}},\\\\
& c_{n+2}(x)=\displaystyle{\frac{x}{3\cdot 2^{1/3}}}c_{n+1}(x)+\displaystyle{\frac{n+1}{3\cdot 2^{2/3}}}c_n(x),\hskip 2cm n=0,1,2,...
\end{array}
$$
%
It is clear that the integral \Eqref{negli} is exponentially small compared with $P_1(x,y)+P_2(x,y)$ in case 3. In summary, we have that
\begin{equation}\label{ptotal}
P(x,y)=\left\lbrace
\begin{array}{cl}
P_1(x,y)+P_2(x,y)+\mathcal{O}(e^{-1.38077\vert y\vert^{4/3}}) & \hskip 1cm {\rm when}\hskip 2mm \displaystyle{\vert\arg y\vert\le\frac{\pi}{8}}, \\
P_1(x,y) & \hskip 1cm {\rm when}\hskip 2mm \displaystyle{-\frac{\pi}{2}\le\arg y<-\frac{\pi}{8}}, \\
P_2(x,y) & \hskip 1cm {\rm when}\hskip 2mm \displaystyle{\frac{\pi}{8}<\arg y\le\frac{\pi}{2}}.
\end{array}
\right.
\end{equation}
where $P_k(x,y)$; $k=1,2$, have the asymptotic expansion \Eqref{puno}, with coefficients given in \Eqref{pdos}, \Eqref{ce} and \Eqref{aa}.

\section{Final remarks and numerical experiments}

For $\Im y>0$, the real part of the dominant term in the exponential factor of $P_2(x,y)$ is larger than the real part of the dominant term in the exponential factor of  $P_1(x,y)$. Therefore, indeed, only $P_2(x,y)$ contributes to the asymptotic behavior of the Pearcy integral: $P(x,y)\sim P_2(x,y)$ for $\Im y>0$. When $\Im y<0$ the situation is just the opposite one, and only $P_1(x,y)$ contributes to the asymptotic behavior of the Pearcy integral:  $P(x,y)\sim P_1(x,y)$ for $\Im y<0$. For positive $y$, the real parts of the dominant terms in the exponential factors of  $P_1(x,y)$ and $P_2(x,y)$ are equal and both terms contribute with equal intensity to the asymptotic behavior of the Pearcy integral: $P(x,y)\sim P_1(x,y)+P_2(x,y)$ for $y>0$.

Then, in the complex $y$ plane, the line $\arg y=0$ is the Stokes line, where both terms $P_1(x,y)$ and $P_2(x,y)$ are equally significant. On the other hand, the difference between the the real part of the dominant terms in the exponential factors of $P_1(x,y)$ and $P_2(x,y)$ is maximal at $\arg y=-3\pi/8$ and minimal at $\arg y=3\pi/8$. This means that the anti-Stokes lines are the ray $\arg y=-3\pi/8$, where $P_1(x,y)$ maximally dominates $P_2(x,y)$ and the ray $\arg y=3\pi/8$, where $P_2(x,y)$ maximally dominates $P_1(x,y)$. This discussion is summarized in Figure 4.

\begin{figure}
\begin{center}
\includegraphics[width=0.5\textwidth]{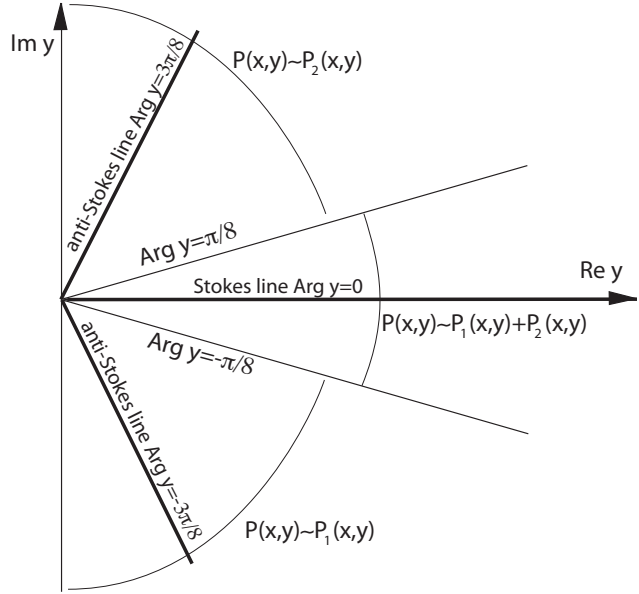}
\end{center}
\caption{\small Stokes and anti-Stokes lines of $P(x,y)$ for $\Re y\ge 0$ and the three different regions for the approximation, as it indicated in \Eqref{ptotal}.}
\label{fig:figura4}
\end{figure}

Finally, we illustrate the accuracy of the expansions \Eqref{puno}-\Eqref{ptotal}. In the following tables we show the relative error for several values of $(x,y)$ and different orders $n$ of the approximation. As we do not have at our disposal the exact value of the Pearcey integral, we have taken the numerical integration of \Eqref{pearcy} obtained with the program {\it Mathematica} with double precision as the exact value of $P(x,y)$.

\begin{table}
\begin{center}
\begin{tabular}{lllllll}\hline
\hskip 3mm$y$ & \multicolumn{6}{c}{$n$} \\
& 0 & 1 & 2 & 3 & 4 & 5\\ \hline
5  & 0.222317& 0.101075& 0.0000918203& 0.00372178& 0.000876593& 0.00302324 \\
10  &0.0316421& 0.00261898& 0.00112219& 0.000403251& 0.0000783942& 0.0000639694 \\
20$e^{i\pi/4}$  &0.0292638 &0.00517274 &0.000228056 &0.0000486543 & 0.0000154281 &
4.73317$\cdot $e-6 \\
20$e^{-3i\pi/8}$  &0.0296318 & 0.00517473 &0.000223576 & 0.0000434364 &0.0000166979 &
5.11767$\cdot$e-6 \\
30  &0.00299077 & 0.00224863 &0.0000906066 &8.36063$\cdot$e-6 &2.84074$\cdot$e-6 &
5.29933$\cdot $e-7 \\
40  &0.0413675 &0.00287761 & 0.0000658777 & 0.0000213951 & 1.41449$\cdot$e-6 &
3.58447$\cdot $e-7 \\
50  &0.0291708 & 0.00152467 & 0.0000388369 &0.0000100058 & 4.79637$\cdot$e-7 &
1.23074$\cdot $e-7 \\
\hline
\end{tabular}
\end{center}
\caption{Relative error for $x=1$ and several values of $y$ and the number of terms $n$ of the approximation given by \Eqref{puno} and \Eqref{ptotal}.}
\label{tab:er3}
\end{table}

\begin{table}
\begin{center}
\begin{tabular}{lllllll}\hline
\hskip 3mm$y$ & \multicolumn{6}{c}{$n$} \\
& 0 & 1 & 2 & 3 & 4 & 5\\ \hline
5  & 0.137947&0.0410408&0.0115823&0.00357474& 0.0159012&0.00749881 \\
10  &0.0443761& 0.0102376& 0.00254929&0.000330121& 0.00115553&0.000371235 \\
20  &0.0312556&0.00192754&0.00045748&0.000173494&0.00006291&
0.0000168014 \\
30$e^{i\pi/4}$ &0.0237833& 0.00108374&0.000209653&0.0000599538& 0.0000165985&
3.36842$\cdot$e-6 \\
30$e^{-3i\pi/8}$ &0.023678& 0.00109324& 0.000206658&0.0000588055&0.0000164115&
3.3194$\cdot$e-6 \\
40  &0.023888&0.00075983& 0.000110305& 0.0000383009&2.30021$\cdot$e-6&
1.02546$\cdot$e-6 \\
50  &0.00206123&0.000599664& 0.0000920699& 6.27624$\cdot$e-7& 3.52299$\cdot$e-6&
5.09017$\cdot$e-7 \\
\hline
\end{tabular}
\end{center}
\caption{Relative error for $x=-2$ and several values of $y$ and the number of terms $n$ of the approximation given by \Eqref{puno} and \Eqref{ptotal}.}
\label{tab:er3}
\end{table}
\medpar



%
\section{Acknowledgments}

The {\it Universidad P\'ublica de Navarra} is acknowledged by its financial support.

\footnotesize{
}
\end{document}